\documentclass[12pt,reqno,british]{amsart}

\usepackage[latin1]{inputenc}
\usepackage[T1]{fontenc}
\usepackage{babel,eucal,url,amssymb,bm,stmaryrd,booktabs}

\theoremstyle{plain}
\newtheorem{theorem}{Theorem}[section]
\newtheorem{proposition}[theorem]{Proposition}

\theoremstyle{definition}

\theoremstyle{remark}
\newtheorem{remark}[theorem]{Remark}
\newtheorem*{acknowledgements}{Acknowledgements}

\numberwithin{equation}{section}

\newcommand{\Lie}[1]{\operatorname{\textsl{#1}}}
\newcommand{\lie}[1]{\operatorname{\mathfrak{#1}}}
\newcommand{\GL}{\Lie{GL}}

\newcommand{\SO}{\Lie{SO}}

\newcommand{\SP}{\Lie{Sp}}
\newcommand{\sP}{\lie{sp}}
\newcommand{\Spin}{\Lie{Spin}}
\newcommand{\SU}{\Lie{SU}}
\newcommand{\su}{\lie{su}}
\newcommand{\Un}{\Lie{U}}
\newcommand{\un}{\lie{u}}

\DeclareMathOperator{\diag}{diag}
\DeclareMathOperator{\im}{Im}
\DeclareMathOperator{\Int}{Int}
\DeclareMathOperator{\re}{Re}
\DeclareMathOperator{\sgn}{sgn}
\DeclareMathOperator{\stab}{stab}
\DeclareMathOperator{\vol}{vol}

\DeclareMathOperator{\RP}{\mathbb RP}
\DeclareMathOperator{\CP}{\mathbb CP}
\DeclareMathOperator{\HP}{\mathbb HP}

\newcommand{\Hodge}{\mathord{\mkern1mu *}}
\newcommand{\hook}{\mathbin{\lrcorner}}

\newcommand{\abs}[1]{\left\lvert #1\right\rvert}
\newcommand{\norm}[1]{\left\lVert #1\right\rVert}
\newcommand{\real}[1]{\left\llbracket #1 \right\rrbracket}
\newcommand{\Span}[1]{\left< #1 \right>}

\newcommand{\coh}[3]{[\,#1\mid #2\mid #3\,]}
\newcommand{\och}[2]{[\,#1\mid #2\,)}
\newcommand{\NoG}{\multicolumn{3}{c}{\textit{No \( \Lie G_2 \) structure}}}

\newenvironment{spmatrix}{\left(\smallmatrix}{\endsmallmatrix\right)}
\newcommand{\subarabic}{\renewcommand{\theequation}{\theparentequation:\arabic{equation}}}

\newcommand{\lb}[1]{\raisebox{-6pt}[12pt][0pt]{#1}}

\begin{document}
\title{Cohomogeneity-one $\Lie G_2$-structures}

\author{Richard Cleyton}
\address[Cleyton]{Department of Mathematics and Computer Science\\
University of Southern Denmark\\
Campusvej 55\\
DK-5230 Odense M\\
Denmark}
\email{cleyton@imada.sdu.dk}

\author{Andrew Swann}
\address[Swann]{Department of Mathematics and Computer Science\\
University of Southern Denmark\\
Campusvej 55\\
DK-5230 Odense M\\
Denmark}
\email{swann@imada.sdu.dk}

\subjclass{(2000) Primary 53C25; secondary 57M50, 57S15}
\keywords{G2, holonomy, weak holonomy, cohomogeneity-one}

\begin{abstract}
  \( \Lie G_2 \)-manifolds with a cohomogeneity-one action of a compact Lie
  group~\( G \) are studied.  For \( G \) simple, all solutions with
  holonomy \( \Lie G_2 \) and weak holonomy \( \Lie G_2 \) are classified.
  The holonomy \( \Lie G_2 \) solutions are necessarily Ricci-flat and
  there is a one-parameter family with \( \SU(3) \)-symmetry.  The weak
  holonomy \( \Lie G_2 \) solutions are Einstein of positive scalar
  curvature and are uniquely determined by the simple symmetry group.
  During the proof the equations for \( \Lie G_2 \)-symplectic and \( \Lie
  G_2 \)-cosymplectic structures are studied and the topological types of
  the manifolds admitting such structures are determined.  New examples of
  compact \( \Lie G_2 \)-cosymplectic manifolds and complete \( \Lie G_2
  \)-symplectic structures are found.
\end{abstract}
\maketitle

\section{Introduction}

A \( \Lie G_2 \)-structure on a seven-dimensional manifold~\( M \) is an
identification of the tangent space with the imaginary octonians.
Equivalently, the geometry is determined by a three-form~\( \phi \) which
at each point is of `generic type', in that it lies in a particular open
orbit for the action of \( \GL(7,\mathbb R) \) (such forms are `stable' in
Hitchin's terminology~\cite{Hitchin:forms}).  The three-form \( \phi
\)~determines a Riemannian metric~\( g \) and hence a Hodge-star
operator~\( \Hodge \).

If \( \phi \) and the four-form \( \Hodge\phi \) are both closed, then \( g
\)~is Ricci-flat and has holonomy contained in~\( \Lie G_2 \).  This is one
of the two exceptional holonomy groups in the Berger classification
(see~\cite{Besse:Einstein,Bryant:status}). The first non-trivial complete
examples were constructed by Bryant \& Salamon
\cite{Bryant-Salamon:exceptional} and compact examples have since been
found by Joyce (first in \cite{Joyce:G2-1,Joyce:G2-2} and more recently
in~\cite{Joyce:holonomy}) and by Kovalev~\cite{Kovalev:G2}.

If \( d \phi = \lambda \Hodge\phi \), for some non-zero constant~\( \lambda
\), then \( g \)~is an Einstein metric of positive scalar curvature and \(
M \)~is said to have weak holonomy~\( \Lie G_2 \).  This terminology was
first introduced by Gray~\cite{Gray:weak}.  Many homogeneous examples are
known.  For example, each Aloff-Walach space \( \SU(3)/\Un(1)_{k,\ell} \),
where \( \Un(1)_{k,\ell} =
\{\diag(e^{ik\theta},e^{i\ell\theta},e^{-i(k+\ell)\theta}) \)\} and \(
k,\ell\) are integers, carries two such metrics (see
\cite{Cabrera-Monar-Swann:G2}).  As \( k \) and \( \ell \) vary, this
family includes infinitely many different homeomorphisms types.  A
classification of the compact homogeneous manifolds with weak holonomy~\(
\Lie G_2 \) is given in~\cite{Friedrich-KMS:G2}.

In this paper we study \( \Lie G_2 \)-structures with a cohomogeneity-one
action of a compact Lie group~\( G \).  This means that \( G \) acts on~\(
M \) preserving the three-form~\( \phi \) and that the generic orbit on~\(
M \) has dimension~\( 7-1=6 \).  We will first determine the connected
groups \( G \) that can act.  Thereafter we study the equations for
holonomy and weak holonomy \( \Lie G_2 \) structures in the case that \( G
\) is simple and determine all solutions.  The simple groups in question
are \( \Lie G_2 \), \( \SP(2) \) and \( \SU(3) \).  In each case we find
that the weak holonomy \( \Lie G_2 \) solutions are unique; they are only
complete in the case with symmetry \( \Lie G_2 \), and here one gets the
round metric on the seven-sphere (and its quotient \( \RP(7) \)). The
limited number of solutions is in strong contrast to the homogeneous case.
For holonomy~\( \Lie G_2 \), the solutions for the first two symmetry
groups are isolated, whereas for \( \SU(3) \) there is a one-parameter
family of solutions.  This family contains a unique complete metric, which
turns out to have \( \Un(3) \)-symmetry.  The \( \Lie G_2 \)-symmetric
solution is flat, whereas those with symmetry \( \SP(2) \) and \( \Un(3) \)
are the metrics found by Bryant \& Salamon
\cite{Bryant-Salamon:exceptional}.  In private communications, Andrew
Dancer \& McKenzie Wang and Gary Gibbons \& Chris Pope tell us that they
have also recently found the one-dimensional family of triaxial \( \SU(3)
\)-symmetric metrics.  Note that by considering non-simple symmetry groups
new complete metrics with holonomy \( \Lie G_2 \) have been found by
Brandhuber et al.~\cite{Brandhuber-GGG:G2}.

Both weak holonomy and holonomy structures satisfy \( d\Hodge\phi=0 \) and
so are special examples of cosymplectic \( \Lie G_2 \)-structures.  Any
hypersurface in an eight-manifold with holonomy \( \Spin(7) \) carries a
cosymplectic \( \Lie G_2 \)-structure and homogeneous cosymplectic \( \Lie
G_2 \)-structures with symmetry \( \SP(2) \) are behind the new \( \Spin(7)
\)-holonomy examples constructed in~\cite{Cvetic-GLP:Spin7}.  Our approach
gives examples of compact cohomogeneity-one manifolds with cosymplectic \(
\Lie G_2 \)-structures.  By Hitchin~\cite{Hitchin:forms} these are
hypersurfaces in manifolds of holonomy \( \Spin(7) \).  It is therefore an
interesting question for future work, which of these \( \Spin(7) \) metrics
are complete.

The other part of the holonomy \( \Lie G_2 \)-equations is \( d\phi=0 \).
Solutions to this equation define what are known as symplectic \( \Lie G_2
\)-structures.  We show that for cohomogeneity-one manifolds with simple
symmetry group, a symplectic \( \Lie G_2 \)-structure exist only if the
manifold also admits a holonomy \( \Lie G_2 \) metric.

\begin{acknowledgements}
  This paper is based on part of the Ph.D. thesis~\cite{Cleyton:thesis} of
  the first named author written under the supervision of the second named
  author.  We thank Andrew Dancer and Anna Fino for useful comments and
  remarks.  Both authors are members of the \textsc{Edge}, Research
  Training Network
  \textsc{hprn-ct-\oldstylenums{2000}-\oldstylenums{00101}}, supported by
  The European Human Potential Programme.
\end{acknowledgements}

\section{$\Lie G_2$-Structures}

Let \( W \) be \( \mathbb R^7 \) with it usual inner product~\( g_0 \).
Take \( \{v^0,\dots,v^6\} \) to be an orthonormal basis for~\( W \) and
write \( v_{01} = v_0v_1 = v_0\wedge v_1 \), etc., in the exterior
algebra~\( \Lambda^*W^* \).  For each \( \theta\in\mathbb R \), we define a
three-form~\( \phi(\theta) \) on~\( W \) by
\begin{equation}
  \label{eq:phi-theta}
  \phi(\theta) = \omega_0 \wedge v_0 + \cos\theta\,\alpha_0
  + \sin\theta\,\beta_0,
\end{equation}
where \( \alpha_0 = v_{246}-v_{235}-v_{145}-v_{136} \), \( \beta_0 =
v_{135}-v_{146}-v_{236}-v_{245} \) and \( \omega_0 = v_{12}+v_{34}+v_{56}
\).

The Lie group \( \Lie G_2 \) may be defined to be the stabiliser of \(
\phi(0) \) under the action of~\( \GL(7,\mathbb R) \).  From this Bryant
shows that \( \Lie G_2 \) is a compact, connected, simply-connected Lie
group of dimension~\( 14 \) \cite{Bryant:exceptional}.  The subgroup of~\(
\Lie G_2 \) fixing \( v^0 \) is isomorphic to~\( \SU(3) \).  Indeed in the
basis \( u_0 = v_0 \), \( u_k = v_{2k-1}+ i v_{2k} \), \( k=1,2,3 \), for
\( W^*\otimes\mathbb C \) we have
\begin{equation*}
  \phi(\theta)= \tfrac i2\bigl((u_1\bar u_1+u_2\bar u_2+u_3\bar u_3)u_0
  + e^{-i\theta}u_1u_2u_3- e^{i\theta}\bar u_1\bar u_2\bar u_3\bigr).
\end{equation*}
Thus \( \phi(\theta) = e^{-i\theta/3}\phi(0) \) showing that stabilisers
of~\( \phi(\theta) \) are all conjugate in~\( \SO(7) \) and that \(
6g_0(v,w) \vol_0 = \bigl(v\hook\phi(\theta)\bigr)\wedge
\bigl(w\hook\phi(\theta)\bigr)\wedge\phi(\theta) \) is independent of~\(
\theta \).

Conversely, the Lie group~\( \Lie G_2 \) acts transitively on the unit
sphere in~\( \mathbb R^7 \).  A choice of unit vector \( v^0 \), determines
a stabiliser isomorphic to \( \SU(3) \) and the action of \( \SU(3) \)
on~\( \Span{v^0}^\bot \) fixes a Kähler form \( \omega_0 \) and a complex
volume which may be written as \( e^{i\theta}u_1u_2u_3 \).  In this way, we
see that there is an orthonormal basis so that the \( \Lie G_2 \)
three-form is \( \phi(\theta) \) as in~\eqref{eq:phi-theta}.

A \( \Lie G_2 \)-structure on a seven-dimensional manifold~\( M \) is
specified by fixing a three-form~\( \phi \) such that for each~\( p \)
there is a basis of \( W=T_pM \) so that \( \phi_p = \phi(\theta) \) for
some~\( \theta \).  We say that a compact Lie group~\( G \) acts on \(
(M^7,\phi) \) with cohomogeneity-one if the \( G \) preserves the
three-form~\( \phi \) and the largest \( G \)-orbits are of dimension~\( 6
\).  In this case, \( B=M/G \) is a one-dimensional manifold, quite
possibly with boundary.  The orbits lying over the interior of~\( B \) are
all isomorphic to~\( G/K \), where \( K=K_p \) is the stabiliser of a~\( p
\in M \) with \( G\cdot p\in \Int B \).  We call these orbits
\emph{principal} and any remaining orbits are called \emph{special}.  Let
\( G/H \) be a special orbit.  Using the action of~\( G \) we may assume
that \( H \) is a subgroup of~\( K \).  A necessary and sufficient for \( M
\) to be a smooth manifold is that for each special orbit \( G/H \), the
quotient \( H/K \)~is a sphere~\cite{Mostert}.

\section{Principal Orbit Structure}

The requirement that~\( G \) acts on~\( M \) with cohomogeneity-one
preserving~\( \phi \) implies that the representation of the isotropy group
\( K=K_p \) on the tangent space of a principal orbit is as a subgroup of
\( \SU(3) \) on its standard six-dimensional representation~\(
\real{\Lambda^{1,0}} \cong \mathbb R^6 \).  Considering the Lie algebras
only we find that \( \lie k \) must be isomorphic to either \( \su(3) \),
\( \un(2) \), \( \su(2) \), \( 2\un(1) \), \( \un(1) \) or~\( \{0\} \).
The possible isotropy representations are then the real representations
underlying the following three-dimensional complex representations: the
standard representation of~\( \su(3) \), the representation \(
L^2\oplus\bar{L} V \) of~\( \un(2) \), the representations \( S^2V \) and
\( \mathbb C\oplus V \) of~\( \su(2) \), the representation \( L_1\oplus
L_2\oplus\bar L_1\bar L_2 \) of~\( 2\un(1) \), the representation \(
L\oplus\bar L \oplus\mathbb C \) of~\( \un(1) \) and finally the trivial
representation \( 3\mathbb C \) of~\( \{0\} \). For each of the non-trivial
representations~\( U \) of a possible isotropy algebra \( \lie k \) the
direct sum~\( \lie g=\lie k\oplus U \) happens to determine a unique
compact real Lie algebra.  These are, respectively, \( \lie g_2 \), \(
\sP(2) \), \( 3\su(2) \), \( \su(3)\oplus\un(1) \), \( \su(3) \) and \(
2\su(2)\oplus\un(1) \).  The trivial representation may be taken to
represent either \( 2\su(2) \), \( \su(2)\oplus 3\un(1) \) or~\( 6\un(1)
\).

If, on the other hand, \( G/K \) is any effective \( 6 \)-dimensional
homogeneous space with \( K \) acting on the isotropy representation as a
subgroup of~\( \SU(3) \), then we may pick an invariant Kähler form \(
\omega \) and an invariant complex volume form \( \alpha \) on \( G/K \)
and obtain a non-degenerate three-form on \( M=\mathbb R\times G/K \) by
defining \( \phi=dt\wedge \omega+\re(\alpha) \).

\begin{theorem}
  Let \( (M^7,\phi) \) be a \( \Lie G_2 \)-manifold of cohomogeneity-one
  under a compact, connected Lie group.  Then, as almost effective
  homogeneous spaces, the principal orbits are one of the following:
  \begin{gather*}
    S^6=\frac{\Lie G_2}{\SU(3)},\qquad
    \CP(3) = \frac{\SP(2)}{\SU(2)\Un(1)},\qquad
    F_{1,2}=\frac{\SU(3)}{T^2},  \\
    S^3\times S^3=\frac{\SU(2)^3}{\SU(2)} =
    \frac{\SU(2)^2 T^1}{T^1} =
    \SU(2)^2, \\
    S^5\times S^1=\frac{\SU(3) T^1}{\SU(2)},\qquad
    S^3\times (S^1)^3=\SU(2) T^3,\qquad
    (S^1)^6=T^6,
  \end{gather*}
  up to finite quotients.  Conversely, any cohomogeneity-one manifold with
  one of these as principal orbit carries a \( \Lie G_2 \)-structure.
  \qed
\end{theorem}

In this paper we will consider the case when \( G \) is simple.  The
principal orbits are the first three cases listed above.  The first of
these is distinguished from the other two in that \( K \) acts irreducibly
on~\( U \).

\section{Irreducible Isotropy}
\label{sec:G2-sym}
This is the case when the principal orbit is~\( \Lie G_2/\SU(3) \).  The
isotropy representation is the real module underlying the standard
representation~\( \Lambda^{1,0} \cong \mathbb C^3 \) of~\( \SU(3) \).  Up
to scale this admits precisely one invariant two-form~\( \omega \) and one
invariant symmetric two-tensor~\( g_0 \).  The space of invariant
three-forms is two-dimensional, spanned by \( \alpha \) and~\( \beta \).
We fix the scales as follows.  Set \( g_0 \) to be the canonical metric
on~\( S^6 = \Lie G_2/\SU(3) \) with sectional curvature one.  Then let \(
\omega \), \( \alpha \) and~\( \beta \) be such that \( \omega^3=6\vol_0
\), \( d\omega=3\alpha \), \( \Hodge_0\alpha=\beta \) and \(
d\beta=-2\omega^2 \).

Let \( \gamma \) be a geodesic through~\( p \) orthogonal to the principal
orbit \( \Lie G_2/\SU(3) \) and parameterize~\( \gamma \) by arc length~\(
t \in I \subset \mathbb R \).  Then the union of principal orbits is \(
I\times \Lie G_2/\SU(3) \subset M \) and there are smooth functions \(
f,\theta \colon I \to \mathbb R \) such that
\begin{gather}
  \label{eq:G2-g}
  g = dt^2+f^2 g_0,\qquad \vol = f^6\vol_0 \wedge dt,\\
  \label{eq:G2-phi}
  \phi=f^2\omega\wedge dt +f^3(\cos\theta\,\alpha+\sin\theta\,\beta).
\end{gather}
Note that \( f(t) \)~is non-zero for each~\( t\in I \).  Our choice of
scales now gives
\begin{gather*}
  \Hodge\phi = \tfrac 12 f^4\omega^2 +
  f^3(\cos\theta\,\beta-\sin\theta\,\alpha) \wedge dt,\\
  d\Hodge\phi = 2f^3\left(f' - 4f^3\cos\theta\right)\omega^2 \wedge dt,\\
  d\phi = \left(3f^2-(f^3\cos\theta)'\right)\alpha \wedge dt
  - (f^3\sin\theta)'\beta \wedge dt
  - 2f^3\sin\theta\,\omega^2.
\end{gather*}

We first consider the cosymplectic \( \Lie G_2 \)-equations \(
d\Hodge\phi=0 \) which are equivalent to \( f'=\cos\theta \).  Locally,
these are described by the one arbitrary function~\( \theta \).
Alternatively, one may regard them as determined by solutions to the
differential inequality \( \abs{f'}\leqslant 1 \).

Geometrically the solutions may be understood as follows.  Consider \(
\mathbb R^8 = W\times \mathbb R \) with its standard \( \Spin(7) \)
four-form \( \Omega = \phi(0)\wedge v_8 + \Hodge_7 \phi(0) \).  As \(
\Spin(7) = \stab_{\GL(8,\mathbb R)}\Omega \) acts transitively on the unit
sphere in~\( \mathbb R^8 \) with stabiliser~\( \Lie G_2 \), we see that for
any unit vector~\( N \), the three-form \( N\hook\Omega \) defines a \(
\Lie G_2 \)-structure on \( \Span{N}^\bot \) and that \( \Omega = \phi
\wedge N^\flat + \Hodge \phi \).  As \( \Omega \)~is closed we therefore
have Gray's observation that any oriented hypersurface~\( H \subset \mathbb
R^8 \) with unit normal~\( N \) carries a cosymplectic \( \Lie G_2
\)-structure.

The hypersurface \( H = \{ (v,s) \in W\times \mathbb R : \norm v = r(s) \}
\) is of cohomogeneity one under the action of~\( \Lie G_2 \).  Its metric
is~\( (1+ (dr/ds)^2) ds^2 + r^2 g_0 \).  Reparameterizing so that \( dt =
\sqrt{(1+(dr/ds)^2)} \,ds \), we obtain a metric in the
form~\eqref{eq:G2-g} with \( f(t)=r(s(t)) \) and hence \( f'(t) =
(dr/ds)/\sqrt{(1+(dr/ds)^2)} \).  However, this has \( \abs{f'(t)}<1 \), so
we may write \( f'=\cos\theta \) and we see that locally each cosymplectic
\( \Lie G_2 \)-solution is given this way away from \( \abs{\cos\theta}=1
\).

The symplectic \( \Lie G_2 \)-equations \( d\phi=0 \) imply first that \(
\sin\theta\equiv0 \).  We then get \( \abs{\cos\theta}=1 \) and \(
f'=\cos\theta \), so such metrics are also cosymplectic and have
holonomy~\( \Lie G_2 \).  However, the solutions are simply \( f(t)=\pm t
\) and we get the standard flat metric on~\( \mathbb R^7 \) with its
standard \( \Lie G_2 \)-structure.

The equations \( d\phi=\lambda\Hodge\phi \) for weak holonomy \( \Lie G_2
\) give
\begin{equation*}
  \lambda f=-4\sin\theta\qquad\text{and}\qquad 4\theta'=-\lambda.
\end{equation*}
Thus \( f(t)=\tfrac4\lambda\sin\left(\lambda t/4\right) \).  The
hypersurface discussion above shows that this is locally the round metric
on \( S^7 \).

\section{Reducible Isotropy: The Equations}

Let us begin with the case of \( \SU(3) \)-symmetry.  The principal
isotropy group \( K=T^2=S^1_1\times S^1_2 \) acts on the standard
representation \( \Lambda^{1,0} \cong\mathbb C^3 \) as \( L_1+L_2+\bar
L_1\bar L_2 \), where \( L_i\cong\mathbb C \), are the standard
representations of \( S^1_i\cong\Un(1) \).  Using the isomorphism \(
\su(3)\otimes\mathbb C\cong\Lambda^{1,1}_0 \), we find that the isotropy
representation is~\( \real{L_1\bar L_2} + \real{L_1L^2_2} + \real{L_1^2L_2}
\).  Each irreducible submodule carries an invariant metric~\( g_i \) and
symplectic form~\( \omega_i \), \( i=1,2,3 \), but the space of invariant
three-forms has dimension~\( 2 \).  Identifying \( T^2 \) with the diagonal
matrices in~\( \SU(3) \), we fix the basis
\begin{gather*}
  E_1=\tfrac12
  \begin{spmatrix}
    0&0&0\\0&0&-1\\0&1&0
  \end{spmatrix},\quad
  E_3=\tfrac12
  \begin{spmatrix}
    0&0&1\\0&0&0\\-1&0&0
  \end{spmatrix},\quad
  E_5=\tfrac12
  \begin{spmatrix}
    0&-1&0\\1&0&0\\0&0&0
  \end{spmatrix},\\
  E_2=\tfrac 1{2i}
  \begin{spmatrix}
    0&0&0\\0&0&1\\0&1&0
  \end{spmatrix},\quad
  E_4=\tfrac 1{2i}
  \begin{spmatrix}
    0&0&1\\0&0&0\\1&0&0
  \end{spmatrix},\quad
  E_6=\tfrac 1{2i}
  \begin{spmatrix}
    0&1&0\\1&0&0\\0&0&0
  \end{spmatrix}
\end{gather*}
of the tangent space at the origin and let \( \{ e_1,\dots,e_6 \} \) denote
the dual basis.  We may now write
\begin{align*}
  g_1&=e_1^2+e_2^2, & g_2&=e_3^2+e_4^2, & g_3&=e_5^2+e_6^2,\\
  \omega_1&=e_{12}, & \omega_2&=e_{34}, & \omega_3&=e_{56},
\end{align*}
and find that
\begin{equation*}
  \alpha=e_{246}-e_{235}-e_{145}-e_{136},\qquad
  \beta=e_{135}-e_{146}-e_{236}-e_{245}
\end{equation*}
is a basis for the invariant three-forms. Put \( \vol_0 = e_{123456} \).
As left-invariant one-forms on~\( \SU(3) \) we have \(
de_i(E_j,E_k)=e_i([E_j,E_k]) \).  One may thus show that on~\( \SU(3)/T^2
\) one has
\begin{equation}
  \label{eq:d-su3}
  \begin{gathered}
    d\omega_1=d\omega_2=d\omega_3=\tfrac12\alpha,\quad d\alpha=0,\\
    d\beta= -2(\omega_1\omega_2+\omega_2\omega_3+\omega_3\omega_1)
    \quad\text{and}\quad d(\omega_i\omega_j)=0.
\end{gathered}
\end{equation}
Any \( \SU(3) \)-invariant \( \Lie G_2 \)-structure on \( I\times
\SU(3)/T^2 \) has
\begin{equation}
  \label{eq:su3-g-vol}
  g = dt^2 + f_1^2 g_1 + f_2^2 g_2 + f_3^2 g_3,\qquad 
  \vol = f_1^2f_2^2f_3^2 \vol_0\wedge dt,
\end{equation}
where \( t \in I\subset\mathbb R \)~is the arc-length parameter of an
orthogonal geodesic and \( f_i \)~are non-vanishing functions.  Using the
equation \( (X\hook\phi)\wedge(Y\hook\phi)\wedge\phi=6g(X,Y)\vol \) and
normalisation \( \phi\wedge\Hodge\phi=7\vol \) we find that the
corresponding invariant three-form is
\begin{equation}
  \label{eq:su3-phi}
  \phi = (f_1^2\omega_1+f_2^2\omega_2+f_3^2\omega_3)\wedge dt +
  f_1f_2f_3(\cos\theta\,\alpha+\sin\theta\,\beta),
\end{equation}
for some function~\( \theta(t) \).  The \( \Lie G_2 \)-structure now has
\begin{equation*}
  \begin{split}
    \Hodge\phi &= f_2^2 f_3^2\omega_2\omega_3 + f_3^2f_1^2\omega_3\omega_1 +
    f_1^2f_2^2\omega_1\omega_2\\
    &\qquad + f_1 f_2 f_3(\cos\theta\, \beta-\sin\theta\, \alpha)\wedge dt,
\end{split}
\end{equation*}
and hence
\begin{align*}
  d\Hodge\phi
  &= \left((f_2^2 f_3^2)' - 2f_1f_2f_3\cos\theta\right)
  \omega_2\omega_3\wedge dt\\
  &\qquad + \left((f_3^2 f_1^2)' - 2f_1f_2f_3\cos\theta\right)
  \omega_3\omega_1\wedge dt\\ 
  &\qquad + \left((f_1^2 f_2^2)' - 2f_1f_2f_3\cos\theta\right)
  \omega_1\omega_2\wedge dt, \\
  d\phi
  &= \left(\tfrac12(f_1^2 + f_2^2 + f_3^2)- (f_1f_2f_3\cos\theta)'\right)
  \alpha\wedge dt \\
  &\qquad - (f_1 f_2 f_3\sin\theta)'\beta\wedge dt \\
  &\qquad - 2f_1 f_2 f_3\sin\theta\,(\omega_1\omega_2 + 
    \omega_2\omega_3 + \omega_3\omega_1).
\end{align*}
We therefore have that the \( \SU(3) \)-invariant \( \Lie G_2 \)-structure
is cosymplectic if
\begin{equation}
  \label{eq:Su3-cs}
  (f_1^2f_2^2)'=(f_3^2f_1^2)'=(f_2^2f_3^2)'=2f_1f_2f_3\cos\theta. 
\end{equation}
It is \( \Lie G_2 \)-symplectic if 
\begin{equation}
  \label{eq:Su3-s}
  (f_1f_2f_3\cos\theta)' = \tfrac12(f_1^2+f_2^2+f_3^2)
  \quad\text{and}\quad
  f_1f_2f_3\sin\theta=0.
\end{equation}
The equations for weak holonomy \( \Lie G_2 \) are
\begin{subequations}
  \label{eq:Su3-wh}
  \begin{gather}
    \label{eq:Su3-wh-a}
    (f_1f_2f_3\cos\theta)'
    = \tfrac12(f_1^2+f_2^2+f_3^2)+\lambda f_1f_2f_3\sin\theta,\\
    \label{eq:Su3-wh-b}
    (f_1f_2f_3\sin\theta)' = -\lambda f_1f_2f_3\cos\theta,\\
    \label{eq:Su3-wh-c}
    -2f_1f_2f_3\sin\theta = \lambda f_1^2f_2^2 = \lambda f_2^2f_3^2 =
    \lambda f_3^2f_1^2. 
  \end{gather}
\end{subequations}

Let us now consider the case of \( \SP(2) \)-symmetry.  The principal
isotropy group \( K=\Un(1)\times\SP(1) \) acts on the standard
representation \( E\cong\mathbb C^4 \) as \( E\cong H+L+\overline L \),
where \( H\cong\mathbb C^2 \) and \( L\cong\mathbb C \) are the standard
representations of~\( \SP(1)=\SU(2) \) and \( \Un(1) \) respectively.
Using \( \sP(2)\otimes\mathbb C\cong S^2E \) we find that the isotropy
representation is~\( \real{L^2}+\real{H\overline L\,} \).  Both of these
modules carry an invariant metric~\( g_i \) and symplectic form~\( \omega_i
\).  The space of invariant three-forms on their sum is two-dimensional.
We give the isotropy representation the basis
\begin{gather*}
  E_1=\tfrac12
  \begin{spmatrix}
    0&0\\0&j
  \end{spmatrix},\quad
  E_2=\tfrac12
  \begin{spmatrix}
    0&0\\0&-k
  \end{spmatrix},\quad
  E_3=\tfrac1{2\sqrt2}
  \begin{spmatrix}
    0&-1\\1&0
  \end{spmatrix},\\
  E_4=\tfrac1{2\sqrt2}
  \begin{spmatrix}
    0&i\\i&0
  \end{spmatrix},\quad
  E_5=\tfrac1{2\sqrt2}
  \begin{spmatrix}
    0&j\\j&0
  \end{spmatrix},\quad
  E_6=\tfrac1{2\sqrt2}
  \begin{spmatrix}
    0&k\\k&0
  \end{spmatrix}.
\end{gather*}
Then the dual elements \( \{e_1,\dots,e_6\} \) are such that \( \{e_1,e_2\}
\) is a basis for \( \real{L^2}^* \) and \( \{e_3,\dots,e_6\} \) is a basis
for \( \real{H\overline L\,}^* \).  We scale \( g_i \) and \( \omega_i \)
so that
\begin{align*}
  g_1&=e_1^2+e_2^2,& g_2&=e_3^2+e_4^2+e_5^2+e_6^2,\\
  \omega_1&=e_{12},& \omega_2&=e_{34}+e_{56}.
\end{align*}
Then
\begin{equation*}
  \alpha=e_{246}-e_{235}-e_{145}-e_{136},\qquad
  \beta=e_{135}-e_{146}-e_{236}-e_{245}
\end{equation*}
is a basis for the invariant three-forms.  Put \( \vol_0 = e_{123456} \).
Using the Lie algebra structure of~\( \sP(2) \) one finds that the
corresponding left-invariant forms on \( \SP(2)/(\Un(1)\times\SP(1)) \)
satisfy
\begin{gather*}
  d\omega_1=\tfrac12\alpha,\qquad d\omega_2=\alpha,\\
  d\alpha=0\quad\text{and}\quad d\beta=-2\omega_1\omega_2-\omega_2^2.
\end{gather*}
Proceeding as in the \( \SU(3) \)-case one finds that the \( \SP(2)
\)-invariant \( \Lie G_2 \)-structures are given by
equations~\eqref{eq:su3-g-vol} and~\eqref{eq:su3-phi} with \( f_3 \equiv
f_2 \).  Computing further, one finds that the equations for these
structures to be cosymplectic, symplectic or have weak holonomy~\( \Lie G_2
\) are those for \( \SU(3) \)-symmetry with \( f_3 \equiv f_2 \).  We may
therefore treat \( \SP(2) \)-symmetry as if it were a special case of~\(
\SU(3) \)-symmetry.

\section{Solving the Cosymplectic \( \Lie G_2 \) Equations}

Consider the cosymplectic \( \Lie G_2 \) equations~\eqref{eq:Su3-cs}.  The
differences of the differentials gives that \( f_i^2(f_j^2-f_k^2) \) is
constant for any permutation \( (i j k) \) of~\( (123) \).  We may
therefore relabel the~\( f_i \) so that \( f_3^2\geqslant f_2^2 \geqslant
f_1^2 \geqslant 0 \) for all~\( t \) and write
\begin{equation}
  \label{eq:mu-nu}
  f_1^2(f_3^2-f_2^2) = \mu^2,\quad 
  f_2^2(f_3^2-f_1^2) = \nu^2,\quad 
  f_3^2(f_2^2-f_1^2) = \nu^2 - \mu^2,
\end{equation}
for some constants \( \nu\geqslant\mu\geqslant0 \).

Let us first deal with two special cases.  If \( \nu=0 \), then \(
f_1^2=f_2^2=f_3^2 \) and we are left with the equation
\begin{equation*}
  2f_1'= \pm \cos\theta.
\end{equation*}
Up to a factor of~\( 2 \) this is just the equation obtained for \( \Lie
G_2 \)-symmetry in~\S\ref{sec:G2-sym}.  Note that we have \(
\abs{f_1'}\leqslant1/2 \).

If \( \nu>\mu=0 \), then \( 2f_2^2 = f_1^2 + \sqrt{f_1^4+4\nu^2} \) and \(
f_1'= \cos\theta\,(1+f_1^2/\sqrt{f_1^4+4\nu^2})^{-1} \), with \( \theta
\)~an arbitrary function.  Note that in this case \( \abs{f_1'}\leqslant 1
\) and \( \abs{f_2'}=\abs{f_1\cos\theta/2f_2}<1/2 \).  

The general case is \( \nu\geqslant\mu>0 \).  Here \( f_3^2>f_2^2\geqslant
f_1^2>0 \) and Equations~\eqref{eq:mu-nu} may be rearranged to give
\begin{subequations}
  \label{eq:f-q}
  \subarabic
  \begin{gather}
    \label{eq:f23}
    f_2^2+\nu^2 f_2^{-2} = f_3^2+(\nu^2-\mu^2)f_3^{-2},\\
    \label{eq:f31}
    f_3^2-(\nu^2-\mu^2)f_3^{-2} = f_1^2+\mu^2f_1^{-2},\\
    \label{eq:f12}
    f_1^2-\mu^2f_1^{-2}=f_2^2-\nu^2f_2^{-2}.
  \end{gather}
\end{subequations}
Regarding equations~\eqref{eq:f-q} as quadratic in~\( f_i^2 \), one sees
that the corresponding discriminants are non-negative.

Let \( \Delta(i;j) \) be the discriminant of~(\ref{eq:f-q}:i) with respect
to~\( f_j^2 \).  Then we have 
\begin{align*}
  \Delta_1 :=
  \Delta(2;3)
  &= \left(f_1^2+\mu^2f_1^{-2}\right)^2+4(\nu^2-\mu^2)\\
  &= \left(f_1^2-\mu^2f_1^{-2}\right)^2+4\nu^2
  = \Delta(3;2) \\
  &= \left(f_3^2+(\nu^2-\mu^2)f_3^{-2}\right)
  \left(f_2^2+\nu^2f_2^{-2}\right), 
  \\
  \Delta_2 :=
  \Delta(3;1)
  &= \left(f_2^2-\nu^2f_2^{-2}\right)^2+4\mu^2\\
  &= \left(f_2^2+\nu^2f_2^{-2}\right)^2-4(\nu^2-\mu^2)
  = \Delta(1;3) \\
  &= \left(f_1^2+\mu^2f_1^{-2}\right)
  \left(f_3^2-(\nu^2-\mu^2)f_3^{-2}\right),
  \\
  \Delta_3 :=
  \Delta(1;2)
  &= \left(f_3^2+(\nu^2-\mu^2)f_3^{-2}\right)^2-4\nu^2\\
  &= \left(f_3^2-(\nu^2-\mu^2)f_3^{-2}\right)^2-4\mu^2
  = \Delta(2;1) \\
  &= \left(f_2^2-\nu^2f_2^{-2}\right)
  \left(f_1^2-\mu^2f_1^{-2}\right).
\end{align*}
The positivity of~\( \Delta_3 \) written as \( \Delta(1;2) \) implies that
\( f_3^4-2\nu f_3^2+\nu^2\geqslant\mu^2 \) which in turns gives either \(
f_3^2 \leqslant \nu-\mu \) or \( f_3^2 \geqslant \nu+\mu \).
However, equation~\eqref{eq:f31} implies that \( f_3^4>\nu^2-\mu^2 =
(\nu+\mu)(\nu-\mu)>(\nu-\mu)^2 \), so
\begin{equation*}
  f_3^2 \geqslant \nu+\mu.
\end{equation*}
Also equation~\eqref{eq:f12} implies that \( \varepsilon =
\sgn(f_1^2-\mu) = \sgn(f_2^2-\nu) \) is well-defined.  Using these
remarks we can choose consistent branches of square roots in solving the
quadratic equations~\eqref{eq:f-q}.  For example
solving~\eqref{eq:f12} for~\( f_2^2 \) and writing the discriminant as a
function of~\( f_1^2 \), we get 
\begin{equation*}
  \begin{split}
    (f_1^2f_2^2)'
    &=\tfrac12\left(f_1^4+f_1^2\sqrt{\Delta_1}-\mu^2\right)'\\
    &=2\left(f_1^4+f_1^2\sqrt{\Delta_1}-\mu^2+2\nu^2\right)
    f_1^3f_1'/\sqrt{\Delta_1}\\
    &=4(f_1^2f_2^2+\nu^2)f_1^3f_1'/\sqrt{\Delta_1}
    =4f_2^2f_3^2f_1^3f_1'/\sqrt{\Delta_1}.
  \end{split}
\end{equation*}
Doing similar computations for the other \( (f_i^2f_j^2)' \) and putting
the results in to~\eqref{eq:Su3-cs} gives
\begin{subequations}
  \label{eq:f-prime}
  \subarabic
  \begin{align}
    \label{eq:f1-p}
    f_1'
    &= \tfrac12 f_2^{-1}f_3^{-1}\cos\theta\sqrt{\Delta_1}
    \notag\\
    &= \tfrac12\varepsilon_{23}\cos\theta
    \sqrt{\left(1+\nu^2f_2^{-4}\right)\left(1+(\nu^2-\mu^2){f_3^{-4}}\right)}
    ,\\
    \label{eq:f2-p}
    f_2'
    &= \tfrac12 f_3^{-1}f_1^{-1}\cos\theta\sqrt{\Delta_2}
    \notag\\
    &= \tfrac12\varepsilon_{31}\cos\theta
    \sqrt{\left(1-(\nu^2-\mu^2){f_3^{-4}}\right)\left(1+\mu^2f_1^{-4}\right)}
    ,\\
    \label{eq:f3-p}
    f_3'
    &= \tfrac12\varepsilon f_1^{-1}f_2^{-1}\cos\theta\sqrt{\Delta_3}
    \notag\\
    &= \tfrac12\varepsilon_{12}^*\cos\theta
    \sqrt{\left(1-\mu^2f_1^{-4}\right)\left(1-\nu^2f_2^{-4}\right)}
    ,
  \end{align}
\end{subequations}
where \( \varepsilon_{ij}=\sgn(f_if_j) \) and \(
\varepsilon_{ij}^*=\varepsilon_{ij}\varepsilon \).  We may rewrite the
right-hand side of equation~\eqref{eq:f1-p} so that it only contains \(
\theta \) and \( f_1 \).  Then for a given function~\( \theta \), we get an
implicit differential equation for~\( f_1 \):
\begin{equation}
  \label{eq:cos}
  f_1'= \varepsilon \cos\theta \sqrt{\Xi(f_1,\mu,\nu)},
\end{equation}
where
\begin{equation}
  \label{eq:Xi}
  \Xi(f_1,\mu,\nu) = 
    \frac{f_1^8 + 2(2\nu^2-\mu^2)f_1^4 + \mu^4}{2 f_1^4\left(f_1^4 +
    (2\nu^2-\mu^2) + \sqrt{f_1^8 + 2(2\nu^2-\mu^2)f_1^4 + \mu^4}\right)}.
\end{equation}
Note that this function \( \Xi(f_1,\mu,\nu) \) is positive and decreasing
with
\begin{equation*}
  \lim_{\abs{f_1}\to\infty} \Xi(f_1,\mu,\nu) = \frac14.
\end{equation*}
Alternatively, the structure may be determined by the function~\( f_1 \):

\begin{theorem}
  Consider a cosymplectic \( \Lie G_2 \)-structure preserved by an action
  of~\( \SU(3) \) of cohomogeneity one.  Then the metric is given by
  equation~\eqref{eq:su3-g-vol}.  Arrange the coefficients so that \(
  f_3^2\geqslant f_2^2 \geqslant f_1^2 \).  Then
  \begin{equation}
    \label{eq:fp-bound}
    \abs{f_1'} \leqslant \sqrt{\Xi(f_1,\mu,\nu)},
  \end{equation}
  for some constants \( \nu\geqslant\mu\geqslant0 \).  
  
  Conversely, any smooth function~\( f_1 \) satisfying the differential
  inequality~\eqref{eq:fp-bound} gives a cosymplectic \( \Lie G_2
  \)-structure with \( f_2 \)~determined by equation~\eqref{eq:f12},
  \( f_3 \) by equation~\eqref{eq:f31} and \( \theta \) by \( f_3\cos\theta
  = (f_1f_2)' \).
  \qed
\end{theorem}

Note that by rescaling and reparameterizing we may rid ourselves of one of
the parameters and, for example, when \( \mu\ne0 \) set either \( \mu \),
\( \nu \) or \( \mu + \nu \) equal to one.

The case of \( \SP(2) \)-symmetry is now obtained by setting either \(
\mu=0 \) or \( \mu=\nu \):

\begin{theorem}
  Consider a cosymplectic \( \Lie G_2 \)-structure preserved by an action
  of~\( \SP(2) \) of cohomogeneity one.  Then the metric is given
  by~\eqref{eq:su3-g-vol} with \( f_3=f_2 \).  The difference \(
  f_1^2-f_2^2 \) has constant sign.  If \( f_1^2\leqslant f_2^2 \), then
  \begin{equation*}
    2f_2^2= f_1^2+\sqrt{f_1^4+4\nu^2}\quad\text{and}\quad
    \abs{f_1'}\leqslant\frac{\sqrt{f_1^4+4\nu^2}}{f_1^2+\sqrt{f_1^4+4\nu^2}} 
  \end{equation*}
  for some \( \nu\geqslant0 \).  If \( f_1^2\geqslant f_2^2 \), then
  \begin{equation*}
    2f_1^2= f_2^2+\nu^2f_2^{-2}\quad\text{and}\quad
    \abs{f_2'}\leqslant\frac{\sqrt{f_2^4+4\nu^2}}{2f_2^2} 
  \end{equation*}
  for some \( \nu\geqslant0 \).

  Conversely, any smooth functions \( f_1 \) and \( f_2 \) satisfying the
  above equations determine a cosymplectic \( \Lie G_2 \)-structure.
  \qed
\end{theorem}

Again, we may rescale and reparameterize to obtain \( \nu=0 \) or~\( 1 \).

\section{Topology and Boundary Conditions}

Let us now turn to discussion of the possible topologies of manifolds with
\( \Lie G_2 \)-structure and a compact simple symmetry group~\( G \) acting
with cohomogeneity one.  General references for the cohomogeneity-one
situation may be found in~\cite{Berard-B:Einstein,Bredon:transformation}.

Let \( M \) be a manifold of cohomogeneity-one under \( G \) with principal
isotropy group~\( K \) and base \( B = M/G \).  The possible topologies
for~\( B \) are homeomorphic to either \( \mathbb R \), \( S^1 \), \(
[0,\infty) \) or \( [0,1] \).  In the first case, \( M \)~is homeomorphic
to the product \( \mathbb R \times G/K \) and an invariant tensor \( \tau
\) on~\( M \) is smooth if and only if \( \tau \) is smooth considered as a
function from \( \mathbb R \) to the space of \( K \)-invariant tensors on
the isotropy representation of the principal orbit.

When the base is a circle, the total space \( M \) is homeomorphic to a
quotient
\begin{equation*}
  \mathbb R\times_h G/K,
\end{equation*}
where \( (t,gK) \) is identified with \( (t+1,ghK) \) for some element \(
h\in N_G(K) \), the normaliser of~\( K \) in~\( G \).  Given \( h \) and \(
h' \) in \( N_G(K) \), these determine the same manifold if \( hK=h'K \)
and they determine equivariantly diffeomorphic manifolds if they satisfy \(
fhf^{-1}=h' \) for some \( f\in N_G(K) \).  For the principal orbits in
question this translates into periodicity requirements corresponding to the
different orders of the elements of \( N_G(K)/K \).  An invariant tensor~\(
\tau_t \) must satisfy
\begin{equation*}
  h^*\tau_t=\tau_{t+1}
\end{equation*}
to be well-defined. 

When the base is a half-open interval, the end-point is the image of a
special orbit with isotropy group~\( H \), where \( H/K \) is diffeomorphic
to a sphere~\( S^m\subset V\simeq\mathbb R^{m+1} \) for some representation
\( V \) of~\( H \).  The total space~\( M \) is then diffeomorphic to the
vector bundle
\begin{equation*}
  M\cong G\times_H V \to G/H.
\end{equation*}
We note that if \( x\in S^m \) has isotropy \( K \) and \( h\in H \)
satisfies \( h\cdot x = -x \) then \( h \) defines an element \( hK\in
N_G(K)/K \) of order~\( 2 \).  Conversely, any non-trivial element \( hK \)
of \( N_G(K)/K \) of order~\( 2 \) defines a subgroup \( H\subset G \) with
\( H/K \) a sphere by taking \( H = K\cup hK \).  An invariant tensor~\(
\tau_t \) on \( M \) must now satisfy
\begin{equation*}
  h^*\tau_t=\tau_{-t}
\end{equation*}
if it is smooth.  This requirement is in general only sufficient when \(
H/K\cong\mathbb Z_2 \).  If \( H/K \) has positive dimension, a metric
two-tensor on \( M_0=M\setminus\pi^{-1}(\{0\}) \) extends to a smooth
metric on \( M \) under the following two conditions.  Firstly, the induced
metric \( g_t(H/K) \) on \( (0,\infty)\times H/K\subset M_0 \) should
satisfy
\begin{equation*}
  g_t(H/K) = dt^2+f^2(t)g_0,
\end{equation*}
where \( g_0 \) is the standard metric on the sphere with sectional
curvature one and \( f \)~is an odd function with \( \abs{f'(0)}=1 \).
Secondly, \( g_t(X,X) \) should be positive everywhere for Killing vector
fields induced by elements of \( \lie h^{\perp} \subset \lie g \).  For the
cases we consider, a \( \Lie G_2 \)-structure on~\( M_0 \) defined by a
three-from \( \phi \) extends to a smooth \( \Lie G_2 \)-structure on \( M
\) if and only if \( h^*\phi_t=\phi_{-t} \) and the metric defined by~\(
\phi \) extends to a smooth metric on~\( M \), see \S\ref{sec:smooth}.

Finally, consider the situation where \( B \) is a closed interval.  Let \(
\pi\colon M\to B \) be the projection.  Then the subspaces \( M_0 =
\pi^{-1}[0,1) \) and \( M_1 = \pi^{-1}(0,1] \) are diffeomorphic to vector
bundles \( G \times_{H_i} V_i \to G/H_i \), where \( H_i \) acts
transitively on the unit sphere in \( V_i \) with isotropy~\( K \).  Given
\( G \), \( K \), \( H_0 \) and~\( H_1 \), the possible diffeomorphism
types of~\( M \) with principal isotropy group~\( K \) and special isotropy
groups \( H_0 \) and~\( H_1 \) are parametrized by the double coset space
\( N_0\backslash N_G(K)/N_1 \), where \( N_i := N_G(K) \cap N_G(H_i) \).
These double cosets correspond to the different equivariant identifications
we may make of \( M_0\setminus\pi^{-1}\{0\} \) with \(
M_1\setminus\pi^{-1}\{1\} \).  The boundary conditions on tensors in this
case are obtained from those for the case of one singular orbit by
considering their restrictions to the half-open intervals.

We will employ the following notation for spaces of of cohomogeneity-one
with special orbits.  When the base \( M/G \) is homeomorphic to the
half-open interval we write \( M=\och{G/H}{G/K} \) where \( G/H \) is the
special orbit over the end point and \( G/K \) is the principal orbit.
When the base is a closed interval we write \( M=\coh{G/H_0}{G/K}{G/H_1}
\).

We now turn to more detailed consideration of our particular principal
orbit types.

\section{Solutions: Irreducible Isotropy}
This is the case of symmetry~\( \Lie G_2 \) with principal orbit \( \Lie
G_2/\SU(3) \).  The normaliser of~\( \SU(3) \) is
\begin{gather*}
  N_{\Lie G_2}(\SU(3)) = \SU(3)\bigcup D_7\SU(3)
\end{gather*}
where \( D_7 = \diag{(-1,1,-1,1,-1,1,-1)} \).  To each of the two elements
of \( N_{\Lie G_2}(\SU(3))/\SU(3) \) corresponds a quotient \( \mathbb R
\times_h \Lie G_2/\SU(3) \) with base diffeomorphic to a circle.

There are precisely two special orbit types: \( \RP(6) = \Lie G_2/N_{\Lie
G_2}(\SU(3)) \) and a point \( \{*\} = \Lie G_2/\Lie G_2 \).  To these
correspond firstly two spaces with base homeomorphic to~\( [0,\infty) \).
The first is the canonical line bundle over~\( \RP(6) \), the second is \(
\mathbb R^7 \) viewed as a \( 7 \)-dimensional vector bundle over a point.

There are three spaces with \( B = [0,1] \) corresponding to the three
possible choices of two special orbits.  If both special orbits are points
the space in question is~\( S^7 \); when one is a point and the other is \(
\RP(6) \) the space is diffeomorphic to \( \RP(7) \); and when both are \(
\RP(6) \) we obtain the connected sum \( \RP(7)\#\RP(7) \).  The
corresponding double coset spaces have precisely one element and therefore
there is only one diffeomorphism type in each case.  The action of \( D_7
\) on the invariant tensors of \( S^6 \) is
\begin{equation*}
  D_7^*(g_0,\omega,\alpha,\beta)=(g_0,-\omega,-\alpha,\beta).
\end{equation*}
As a consequence \( D_7^*\vol_0 = -\vol_0 \).  In particular, the space \(
\mathbb R\times_{D_7}S^6 \) is not orientable and therefore cannot carry a
\( \Lie G_2 \)-structure.

When \( M \) has a special orbit with isotropy~\( \Lie G_2 \) at~\( t=0 \)
the metric \( g_t \) extends to a smooth metric on a neighbourhood of the
special orbit if and only if the function~\( f \) is odd with \(
\abs{f'(0)}=1 \).  The requirement \( D_7^*\phi_t=\phi_{-t} \) now implies
that \( \sin\theta \)~is odd and \( \cos\theta \)~is even around \( t=0 \).
If, on the other hand, the special orbit at \( t=0 \) is~\( \RP(6) \), then
\( f \) must be even and non-zero everywhere for the metric to extend
smoothly.  In that case \( \cos\theta \) must be even and \( \sin\theta \)
odd.

Now consider the cosymplectic equations.  One solution is given by \(
f\equiv c \), where \( c \) is a positive constant, and \( \theta\equiv 0
\).  This solution satisfies the boundary conditions for \(
\coh{\RP(6)}{S^6}{\RP(6)} \) and \( \och{\RP(6)}{S^6} \), as well as the
periodicity requirement for \( \mathbb R\times_{e}S^6=S^1\times S^6 \).

The unique solution to the symplectic equation satisfies the boundary
conditions only for \( \mathbb R^7=\och{*}{S^6} \).  The weak holonomy
solutions \( f(t) = 4\lambda^{-1}\sin(\lambda t/4) \) are smooth on \(
S^7=\coh{*}{S^6}{*} \) for \( t\in[0,4\pi/\lambda] \) and on \(
\RP(7)=\coh{*}{S^6}{\RP(6)} \) for \( t\in[0,2\pi/\lambda] \).  Different
choices of \( \lambda \) scale the metric by a homothety.

\begin{theorem}
  Let \( M^7 \) be a manifold with \( \Lie G_2 \)-structure preserved by an
  action of \( \Lie G_2 \) of cohomogeneity one.  Then the principal orbit
  is \( \Lie G_2/\SU(3) \) and \( M^7 \) is listed in Table~\ref{tab:g2}.
  The symplectic~\( \Lie G_2 \), holonomy~\( \Lie G_2 \) and weak
  holonomy~\( \Lie G_2 \) solutions are unique up to scale.  The first two
  are flat, the last has constant curvature.  \qed
\end{theorem}

\begin{table}[tp]
\begin{tabular}{@{}l@{}ccc@{}}
\toprule
\lb{\( M^7 \)}&Holonomy \&&Weak&\lb{Cosymplectic}\\
&Symplectic&Holonomy&\\
\midrule
\( S^7 \)&None&Complete&Complete\\
\( \RP(7) \)&None&Complete&Complete\\
\( \RP(7)\#\RP(7) \)&None&None&Complete\\
\( S^1\times S^6 \)&None&None&Complete\\
\( \mathbb R\times_{D_7} S^6 \)&\NoG\\[1ex]
\( \och{\RP(6)}{S^6} \)&None&None&Complete\\
\( \mathbb R^7 \)&Complete&Incomplete&Complete\\[1ex]
\( \mathbb R\times S^6 \)&Incomplete&Incomplete&Complete\\
\bottomrule
\end{tabular}
\caption{\( \Lie G_2 \) solutions with symmetry~\( \Lie
G_2 \).}
\label{tab:g2}
\end{table}

\section{Solutions: Reducible Isotropy}
We first consider the instance of \( \SU(3) \)-symmetry; that for \( \SP(2)
\) will then follow relatively easily.  The principal orbits are \(
\SU(3)/T^2 \) and the normaliser of \( T^2 \) in \( \SU(3) \) is
\begin{gather*}
  N_{\SU(3)}(T^2)=\bigcup_{\sigma\in \Sigma_3} A_{\sigma}T^2,
\end{gather*}
where \( \Sigma_3 \) is the symmetric group on three elements
and 
\begin{equation*}
  A_{(123)}=
  \begin{spmatrix}
    0&0&1\\1&0&0\\0&1&0
  \end{spmatrix},\quad\text{and}\quad
  A_{(23)}=
  \begin{spmatrix}
    -1&0&0\\
    0&0&1\\
    0&1&0
  \end{spmatrix}.
\end{equation*}
Therefore there are three spaces \( \mathbb R\times_{A_\sigma} F_{1,2} \)
over the circle corresponding to \( A_{(23)} \), \( A_{(123)} \) and~\( e
\).

There are two special orbit types: \( \CP(2)_1 = \SU(3)/\Un(2)_{(23)} \)
and \( \mathcal F_{(23)}=F_{1,2}/A_{(23)} = \SU(3)/(T^2\cup A_{(23)}T^2)
\), where \( \Un(2)_{(23)} \) is the \( \Un(2)\subset\SU(3) \) containing
\( T^2\cup A_{(23)}T^2 \).  Corresponding to these there are two spaces
with base homeomorphic to the half-open interval.

The double coset spaces \( N_0\backslash N_{\SU(3)}(T^2)/N_1 \) all have
two components.  Therefore we have six different cohomogeneity-one spaces
with the closed interval as base.  When both special orbits are complex
projective spaces, we may write the space from the trivial double coset as
\( \coh{\CP(2)_1}{F_{1,2}}{\CP(2)_1} \) and that from the non-trivial
double coset as \( \coh{\CP(2)_1}{F_{1,2}}{\CP(2)_2} \), where \( \CP(2)_2
= \SU(3)/\Un(2)_{(13)} \).  We use similar notation in the other cases.

Now consider in turn the actions of the elements \( A_{(23)} \) and~\(
A_{(123)} \).  The element \( A_{(23)} \) acts with order~\( 2 \) and
transforms the \( \SU(3) \)-invariant tensors of~\( F_{1,2} \) as
\begin{equation*}
  A_{(23)}^*(g_1,g_2,g_3,\omega_1,\omega_2,\omega_3,\alpha,\beta)=
  (g_1,g_3,g_2,-\omega_1,-\omega_3,-\omega_2,-\alpha,\beta).
\end{equation*}
This implies that the manifold \( \mathbb R \times_{A_{(23)}} F_{1,2} \)
cannot carry a \( \Lie G_2 \)-structure.  It also leads to boundary
conditions on the metric and three-form for the two types of special orbit.
For special orbit \( \CP(2)_1 \) these translate into
\begin{equation*}
  \left.
    \begin{gathered}
      \text{\( f_1 \) and \( \sin\theta \) are odd functions},\\
      \text{\( \cos\theta \) is an even function},\\
      f^2_2(t)=f^2_3(-t),\quad\abs{f_1'(0)}=1\quad\text{and}\quad
      f_2(0)\ne 0.
    \end{gathered}
    \right\}
\end{equation*}
Those for the special orbit \( \mathcal F_{(23)} \) are
\begin{equation*}
  \left.
  \begin{gathered}
    \text{\( f_1 \) and  \(\sin\theta \) are even functions},\\
    \text{\( \cos\theta \) is an odd function},\\
    f^2_2(t)=f^2_3(-t),\quad\text{and}\quad f_1(0)\ne0\ne f_2(0).
  \end{gathered}
  \right\}
\end{equation*} 
Note that in both cases the product \( f_2f_3 \) is even.

The action of \( A_{(123)} \) on the invariant tensors of \( F_{1,2} \) is
\begin{equation*}
  A_{(123)}^*(g_1,g_2,g_3,\omega_1,\omega_2,\omega_3,\alpha,\beta)
  =(g_2,g_3,g_1,\omega_2,\omega_3,\omega_1,\alpha,\beta)
\end{equation*}
whence the periodicity conditions on \( \mathbb R\times_{A_{(123)}}F_{1,2}
\) state that \( f_1^2(t)=f_2^2(t+1)=f_3^2(t+2) \).  Note that the tensors
\begin{equation*} 
  g_0=g_1+g_2+g_3,\quad\omega_0=\omega_1+\omega_2+\omega_3,
  \quad\alpha\quad\text{and}\quad\beta
\end{equation*} 
all are invariant under
\begin{equation}
  \label{eq:T123}
  T^{(123)} = \bigcup_{\text{\( \sigma
  \) even}} A_\sigma T^2.  
\end{equation}
Therefore \( \mathcal F_{(123)}=\SU(3)/T^{(123)} \) is a second possible
principal orbit for symmetry~\( \SU(3) \).  It is not hard to check that \(
A_{(123)} \) generates the only possible finite action on the principal
orbits that preserves an \( \SU(3) \)-invariant \( \Lie G_2 \)-structure.
The normaliser of \( T^{(123)} \) in~\( \SU(3) \) is of course \(
T^{(123)}\cup A_{(23)} T^{(123)} \) and \( A_{(23)} \)
acts on the invariant tensors as
\begin{equation*}
  A_{(23)}^*(g_0,\omega_0,\alpha,\beta)=(g_0,-\omega_0,-\alpha,\beta).
\end{equation*}
For the principal orbit \( \mathcal F_{(123)} \), the analysis is now the
same as for the case of \( \Lie G_2 \)-symmetry discussed in the previous
section.

Returning to principal orbit \( F_{1,2} \) we see that taking \(
f_1=f_2=f_3\equiv c \), with \( c \) a positive constant, and \(
\theta\equiv0 \) solves the \( \SU(3) \)-symmetric cosymplectic equations
as well as the periodicity requirement on \( S^1\times F_{1,2} \) and \(
\mathbb R\times_{A_{(123)}} F_{1,2} \) and the boundary conditions on \(
\och{\mathcal F_{(23)}}{F_{1,2}} \), \( \coh{\mathcal
F_{(23)}}{F_{1,2}}{\mathcal F_{(23)}} \) and \( \coh{\mathcal
F_{(23)}}{F_{1,2}}{\mathcal F_{(13)}} \).

Consider the cosymplectic \( \Lie G_2 \)-equations together with the
boundary conditions for either \( \coh{\CP(2)_1}{F_{1,2}}{\CP(2)_2} \) or
\( \coh{\CP(2)_1}{F_{1,2}}{\mathcal F_{(13)}} \).  From~\eqref{eq:mu-nu},
we have that two of the three constants \( \mu \), \( \nu \) and \(
\nu^2-\mu^2 \) must be zero.  But this implies that the third constant is
also zero and that \( f_1^2=f_2^2=f_3^2 \).  The boundary conditions now
give that \( f_1 \)~is both even and odd at \( t=0 \) which clearly can not
be the case.  Thus these spaces do not carry invariant cosymplectic \( \Lie
G_2 \)-structures.

Finally, let us consider the cosymplectic equations on \(
\coh{\CP(2)_1}{F_{1,2}}{\CP(2)_1} \) and \(
\coh{\CP(2)_1}{F_{1,2}}{\mathcal F_{(23)}} \).  Solutions on these spaces
can be obtained by as follows.  Set
\begin{equation*}
  d\theta=(1+\sin^2\theta)^{1/4}dt,
\end{equation*}
and determine the remaining functions via equation~\eqref{eq:cos}.
The metric is then
\begin{equation*}
  g=(1+\sin^2\theta)^{-1/2}
  \left(d\theta^2+\sin^2\theta g_1+(1+\sin^2\theta)(g_2+g_3)\right)
\end{equation*}
and the three-form is
\begin{multline*}
  \phi = (1+\sin^2\theta)^{-3/4}
  \Bigl(
  \bigl( (\sin^2\theta\,\omega_1 + (1+\sin^2\theta)(\omega_2+\omega_3) \bigr)
  d\theta\\ 
  +\sin\theta(1+\sin^2\theta) (\cos\theta\,\alpha+\sin\theta\,\beta) \Bigr).
\end{multline*}
With \( \theta\in [0,\pi] \) these solve the cosymplectic equations and the
boundary conditions for \( \coh{\CP(2)_1}{F_{1,2}}{\CP(2)_1} \).
Restricting \( \theta \) to \( [0,\pi/2] \) we also get a solution on \(
\coh{\CP(2)_1}{F_{1,2}}{\mathcal F_{(23)}} \).

This completes the discussion of the cosymplectic equations under \( \SU(3)
\)-symmetry.  We will return to the holonomy and weak holonomy equations
after discussing the symmetry group \( \SP(2) \).

\begin{theorem}
  Let \( M^7 \) be a manifold with \( \Lie G_2 \)-structure preserved by an
  action of \( \SU(3) \) of cohomogeneity one.  The principal orbit is
  either \( F_{1,2}=\SU(3)/T^2 \) or its \( \mathbb Z_3 \)-quotient \(
  \mathcal F_{(123)} = \SU(3)/T^{(123)} \), see~\eqref{eq:T123}.  The
  possible \( M^7 \) are listed in Table~\ref{tab:su3} together with
  information on the existence of cosymplectic \( \Lie G_2 \)-structures.
  \qed
\end{theorem}

\begin{table}[tp]
\begin{tabular}{@{}l@{}ccc@{}}
\toprule
\lb{\( M^7 \)}&Holonomy \&&Weak&\lb{Cosymplectic}\\
&Symplectic&Holonomy&\\
\midrule
\( \coh{\CP(2)_1}{F_{1,2}}{\CP(2)_1} \)&None&None&Complete\\
\( \coh{\CP(2)_1}{F_{1,2}}{\CP(2)_2} \)&None&None&None\\
\( \coh{\CP(2)_1}{F_{1,2}}{\mathcal F_{(23)}} \)&None&None&Complete\\
\( \coh{\CP(2)_1}{F_{1,2}}{\mathcal F_{(13)}} \)&None&None&None\\
\( \coh{\mathcal F_{(23)}}{F_{1,2}}{\mathcal F_{(23)}} \)&None&None&Complete\\
\( \coh{\mathcal F_{(23)}}{F_{1,2}}{\mathcal F_{(13)}} \)&None&None&Complete\\
\( S^1\times F_{1,2} \)&None&None&Complete\\
\( \mathbb R\times_{A_{(23)}} F_{1,2} \)&\NoG\\
\( \mathbb R\times_{A_{(123)}} F_{1,2} \)&None&None&Complete\\[1ex]
\( \och{\CP(2)_1}{F_{1,2}} \)&Complete&None&Complete\\
\( \och{\mathcal F_{(23)}}{F_{1,2}} \)&None&None&Complete\\[1ex]
\( \mathbb R\times F_{1,2} \)&Incomplete&Incomplete&Complete\\
\midrule
\( \coh{\mathcal F_\Sigma}{\mathcal F_{(123)}}{\mathcal F_\Sigma}
\)&None&None&Complete\\ 
\( \och{\mathcal F_\Sigma}{\mathcal F_{(123)}} \)&None&None&Complete\\[1ex]
\( S^1 \times \mathcal F_{(123)} \)&None&None&Complete\\
\( \mathbb R \times_{A_{(23)}} \mathcal F_{(123)} \)&\NoG\\[1ex]
\( \mathbb R\times \mathcal F_{(123)} \)&Incomplete&Incomplete&Complete\\
\bottomrule
\end{tabular}
\caption{\( G_2 \) solutions with symmetry~\( \SU(3) \).  Here \( \mathcal
F_\sigma = F_{1,2}/A_\sigma \) and \( \mathcal F_\Sigma = F_{1,2}/\Sigma_3
\).}
\label{tab:su3}
\end{table}

The topological analysis in the case of \( \SP(2) \)-symmetry is very
similar to the \( \Lie G_2 \) case for the simple reason that the
normaliser of \( \Un(1)\SP(1) \) again has two components:
\begin{equation*}
  N_{\SP(2)}(\Un(1)\SP(1))=\Un(1)\SP(1)\bigcup D_2\Un(1)\SP(1),
\end{equation*}
where \( D_2= \begin{spmatrix} j&0\\0&1 \end{spmatrix} \).  Therefore we
have two possible spaces \( \mathbb R \times_h \CP(3) \) with base a
circle, and two possible special orbit types:
\begin{gather*}
  S^4=\HP(1)=\frac{\SP(2)}{\SP(1)\times\SP(1)},\\
  \mathcal C=\CP(3)/\mathbb Z_2=\frac{\SP(2)}{N_{\SP(2)}(\Un(1)\SP(1))}.
\end{gather*}

Let us now consider the action of \( D_2 \) on the invariant tensors of \(
\CP(3) \):
\begin{equation*}
  D_2^*(g_1,g_2,\omega_1,\omega_2,\alpha,\beta) =
  (g_1,g_2,-\omega_1,-\omega_2,\alpha,-\beta)
\end{equation*}
This means that the boundary conditions are precisely those for \( \SU(3)
\)-symmetry with \( f_2 \equiv f_3 \).  In particular the compact, complete
solutions to the cosymplectic equations found for \( \SU(3) \)-symmetry
also give solutions for \( \SP(2) \)-symmetry.  The results of the analysis
in this case may be found in Table~\ref{tab:sp2}.  Note that the existence
of an invariant \( \Lie G_2 \)-structure implies that the only possible
principal orbit is~\( \CP(3) \).

\begin{table}[tp]
\begin{tabular}{@{}l@{}ccc@{}}
\toprule
\lb{\( M^7 \)}&Holonomy \&&Weak&\lb{Cosymplectic}\\
&Symplectic&Holonomy&\\
\midrule
\( \coh{S^4}{\CP(3)}{S^4} \)&None&None&Complete\\
\( \coh{S^4}{\CP(3)}{\mathcal C} \)&None&None&Complete\\
\( \coh{\mathcal C}{\CP(3)}{\mathcal C} \)&None&None&Complete\\
\( S^1\times\CP(3) \)&None&None&Complete\\
\( \mathbb R\times_{D_2}\CP(3) \)&\NoG\\[1ex]
\( \och{S^4}{\CP(3)} \)&Complete&None&Complete\\
\( \och{\mathcal C}{\CP(3)} \)&None&None&Complete\\[1ex]
\( \mathbb R\times\CP(3) \)&Incomplete&Incomplete&Complete\\
\bottomrule
\end{tabular}
\caption{\( \Lie G_2 \) solutions with symmetry~\( \SP(2) \).
Here \( \mathcal C \) denotes \( \CP(3)/\mathbb Z_2 \).}
\label{tab:sp2}
\end{table}

Let us now turn to the weak holonomy equations, firstly for \( \SU(3) \).
Equations \eqref{eq:Su3-wh-c} and~\eqref{eq:mu-nu} imply that \( f_1^2 =
f_2^2 = f_3^2 \) and that \( f_i = - \varepsilon_{jk} 2\lambda^{-1}
\sin\theta \).  Equations \eqref{eq:Su3-wh-a} and~\eqref{eq:Su3-wh-b} are
now
\begin{gather*}
  (\theta'+4\lambda)\sin^2\theta\,(4\cos^2\theta-1)=0,\\
  (\theta'+4\lambda)\sin^3\theta\cos\theta=0.
\end{gather*}
As \( f_i \) is non-zero on the principal orbits, we get that \(
\theta'=-4\lambda \).  We deduce that we have the same behaviour for \(
\SP(2) \)-symmetry.

\begin{theorem}
  Up to scale, the spaces \( (0,\pi/2)\times G/K \), with \( G/K=F_{1,2}
  \), \( \mathcal F_{(123)} \) or \( \CP(2) \) admit unique structures with
  weak holonomy~\( \Lie G_2 \) invariant under the action of \( \SU(3) \)
  or~\( \SP(2) \).  The metric and three-forms are
  \begin{gather*}
    g = 4d\theta^2 + \sin^2\theta\, g_0,\\
    \phi = \sin^2\theta\,(\omega_0\wedge d\theta +
    \sin\theta\,(\cos\theta\,\alpha + \sin\theta\,\beta)),
  \end{gather*}
  where \( g_0 = \sum_i g_i \) and \( \omega_0 = \sum_i \omega_i \).  These
  structures are incomplete and do not extend over any special orbits.
  \qed
\end{theorem}

Next we discuss the holonomy solutions.  The second equation
in~\eqref{eq:Su3-s} implies that \( \sin\theta\equiv0 \).  Let \(
\varepsilon_\theta = \cos\theta \).  Then the cosymplectic
equations~\eqref{eq:Su3-cs} show that the first equation
in~\eqref{eq:Su3-s} is automatically satisfied.  We thus have that \( f_1
\) satisfies the differential equation
\begin{equation*}
  f_1' = \varepsilon_\theta \sqrt{\Xi(f_1,\mu,\nu)},
\end{equation*}
where \( \Xi \) is defined in~\eqref{eq:Xi}.  As \(
\Xi(f_1,\mu,\nu)\geqslant 1/4 \), we have that \( \abs{f_1}\geqslant
\tfrac12 t + c \) and so any complete solution has exactly one special
orbit and \( f_1 \)~vanishes on that orbit.

If \( \nu=0 \) then \( f_1^2=f_2^2=f_3^2=\tfrac14 t \), which does not
satisfy any of the boundary conditions for symmetry \( \SU(3) \) or \(
\SP(2) \).

We may now take \( \nu>0 \) and introduce the parameter change \( r(t)^2 =
f_1^2(t)f_3(t)^2 \), with \( r(t)>0 \).  Then \( \abs{r'}=\abs{(f_1f_3)'} =
\abs{f_2} \) is strictly positive.  Using \eqref{eq:Su3-cs} and
\eqref{eq:mu-nu}, we get
\begin{gather*}
  f_1^2 = r\sqrt{({r^2-\mu^2})/({r^2+\nu^2-\mu^2})},\\
  f_2^2 = r^{-1}\sqrt{(r^2-\mu^2)(r^2+\nu^2-\mu^2)},\\
  f_3^3 = r\sqrt{(r^2+\nu^2-\mu^2)/(r^2-\mu^2)}
\end{gather*}
and
\begin{equation*}
  dt^2 = r\,dr^2/\sqrt{(r^2-\mu^2)(r^2+\nu^2-\mu^2)} .
\end{equation*}
These are `triaxial' metrics with holonomy \( \Lie G_2 \) and \( \SU(3)
\)-symmetry.  To be complete there must be a special orbit.  This requires
\( f_1=0 \) and \( f_2,f_3\ne0 \) at \( t=0 \).  The first condition
implies \( r(0)^2=\mu^2 \), the second gives \( \mu=0 \).  Thus this
solution has \( f_2^2(t)=f_3^2(t) \), which is the metric found by Bryant
\& Salamon on the bundle of anti-self-dual two-forms over~\( \CP(2) \)
\cite{Bryant-Salamon:exceptional}.  This solution also defines a structure
with \( \SP(2) \)-symmetry.

\begin{theorem}
  The space \( \mathbb R\times F_{1,2} \) admits a one-parameter family of
  holonomy \( \Lie G_2 \) metrics with \( \SU(3) \)-symmetry.  Only one
  metric extends to a complete metric, and the underlying manifold is \(
  \och{\CP(2)_1}{F_{1,2}} \), the bundle of anti-self-dual two-forms
  over~\( \CP(2) \).
  
  The space \( \mathbb R\times \mathcal F_{(123)} \) admits a unique
  incomplete metric with holonomy~\( \Lie G_2 \) invariant under~\( \SU(3)
  \).

  The space \( \mathbb R\times \CP(3) \) admits two metrics with
  holonomy~\( \Lie G_2 \).  One is incomplete, the other extends to a
  complete metric on~\( \och{S^4}{\CP(3)} \), the bundle of anti-self-dual
  two-forms over~\( S^4 \).
  \qed
\end{theorem}

\begin{remark}
  As Andrew Dancer pointed out to us the substitutions \( dt =
  f_1f_2f_3\,ds \) and \( w_i = f_jf_k \) for each even permutation \(
  (ijk) \) of \( (123) \) reduce the \( \SU(3) \)-symmetric holonomy \(
  \Lie G_2 \)-equations to Euler's equations for a spinning top.  These
  equations may then be solved by elliptic integrals.  However, as this is
  no longer an arc-length parameterization, one now has to work harder to
  determine questions of completeness.
\end{remark}

Finally we consider the equations~\eqref{eq:Su3-s} for a symplectic \( G_2
\) structure with symmetry \( \SU(3) \).  We have \( \sin\theta=0 \).  Put
\( \varepsilon_\theta = \cos\theta \) and  set
\begin{equation*}
  h^3=f_1f_2f_3, \quad x=f_2^{-2}h^2 \quad\text{and}\quad y=f_3^{-2}h^2,
\end{equation*}
so \( x \) and \( y \) are positive.
Equations~\eqref{eq:Su3-s} then give 
\begin{equation}
  \label{eq:hxy}
  6\,\varepsilon_\theta h' = xy + \frac 1x+ \frac 1y.
\end{equation}
On \( (0,\infty)^2 \), the right-hand side has a global minimum at~\( (1,1)
\) and so \( \abs{h'}\geqslant1/2 \).  This implies that there are no
periodic solutions and that any complete solution has exactly one special
orbit.  As \( f_2f_3 \) is even, we also see that \( f_1 \)~vanishes at the
special orbit.  Therefore we have exactly the same topologies as for
holonomy~\( \Lie G_2 \).  Note however that there are more solutions to the
symplectic equations than for holonomy~\( \Lie G_2 \).  A particularly
simple example of this is furnished by setting \( f_1(t)=t \), \(
f_2^2(t)=1+(t/2)^2 \) and \( f_2\equiv f_3 \).  Complete triaxial solutions
may be obtained as follows: begin with the complete \( \Un(3) \)-symmetric
metric with holonomy~\( \Lie G_2 \); hold \( h \)~fixed, make a smooth
deformation of \( x \) on~\( [1,\infty) \) and determine the corresponding
deformation of~\( y \) by~\eqref{eq:hxy}.

\begin{proposition}
  Let \( M^7 \) admit a \( \Lie G_2 \)-structure which is preserved by a
  cohomogeneity-one action of a compact simple Lie group.  Then \( M^7 \)
  admits an invariant symplectic \( \Lie G_2 \)-structure if and only if \(
  M^7 \) admits an invariant metric with holonomy~\( \Lie G_2 \).
  Similarly, complete symplectic \( \Lie G_2 \)-structures only exist on
  manifolds with complete \( \Lie G_2 \) holonomy metrics.  \qed
\end{proposition}

\section{Smoothness of the Three-form}
\label{sec:smooth}
In this section we will briefly indicate how to check that the
three-form~\( \phi \) is smooth once we have \( h^*\phi_t=\phi_{-t} \) and
smoothness of the metric~\( g \).  The only case where significant work is
required is that of special orbit~\( \CP(2) \) under \( \SU(3) \)-symmetry.
The case of special orbit~\( S^4 \) under \( \SP(2) \)-symmetry follows by
similar arguments.

The manifold \( \och{\CP(2)}{F_{1,2}} \) is \( \SU(3) \)-equivariantly
isomorphic to the bundle of anti-self-dual two-forms~\( \Lambda^2_- \)
over~\( \CP(2) \).  Bryant \& Salamon \cite{Bryant-Salamon:exceptional}
showed how to construct holonomy~\( \Lie G_2 \) metrics on~\( \Lambda^2_-
\), but they did not write down the general \( \SU(3) \)-invariant
three-form because they treated all four-manifolds at once.  In the
following, we specialise Bryant \& Salamon's approach in the style
of~\cite{Swann:MathAnn}.

If we write \( \CP(2) = \SU(3)/\Un(2) \), then \( P = \SU(3) \) is a
principal bundle of frames with structure group~\( \Un(2) =
\Un(1)\times_{\mathbb Z/2}\SP(1) \).  Under the action of~\( \Un(2) \), we
have \( \Lambda^{1,0} \cong HL+\bar L^2 \), where \( L\cong\mathbb C \) and
\( H\cong\mathbb C^2 \) are the standard representations of \( \Un(1) \)
and \( \SP(1) \), respectively.  This may be regarded as an identification
not only of representations but also of bundles over~\( \CP(2) \), if to a
representation~\( V \) of~\( \Un(2) \) we associate the bundle, also
denoted~\( V \),
\begin{equation*}
  P \times_{\Un(2)} V,
\end{equation*}
which is \( P\times V \) modulo the action \( (u,\xi)\mapsto(u\cdot
g,g^{-1}\cdot\xi) \).  We then have \( \Lambda^2_-=S^2H\cong\im\mathbb H
\).  Let \( \bm\theta = \theta_0+\theta_1i+\theta_2j+\theta_3k \in
\Omega^1(P,\mathbb H) \) be the canonical one-form.  Write \( \eta=
\eta_1i+\eta_2j+\eta_3k \in \Omega^1(P,\im \mathbb H) \) for the \( \sP(1)
\)-part of the \( \Un(2) \) Levi-Civita connection.  As the Fubini-Study
metric is self-dual and Einstein one finds that
\begin{equation*}
  d\eta + \eta\wedge \eta = c \bar{\bm\theta}\wedge\bm\theta
\end{equation*}
for some positive constant~\( c \) (a positive constant multiple of the
scalar curvature).  If \( x=x_1i+x_2j+x_3k \) is the coordinate on~\(
\im\mathbb H \) then let \( r^2 = x\bar x = -x^2 \).  The one-form
\begin{equation*}
  \bm\alpha = dx + \eta x - x\eta
\end{equation*}
is semi-basic on \( P\times S^2H \).  One may now check that
\begin{gather*}
  \omega_1 = r^{-3}\bm\alpha x \bm\alpha,\\
  \omega_2 = 4c(r^{-1}\bm\theta x\bar{\bm\theta} + \bar{\bm\theta} i
  \bm\theta),\\
  \omega_3 = 4c(r^{-1}\bm\theta x\bar{\bm\theta} - \bar{\bm\theta} i
  \bm\theta),\\
  \alpha = -4c(r^{-1}\bm\theta\bm\alpha\bar{\bm\theta}+r^{-3}\bm\theta
  x\bm\alpha x\bar{\bm\theta}),\\
  \beta = -4cr^{-2}\bm\theta(\bm\alpha x-x\bm\alpha)\bar{\bm\theta}
\end{gather*}
satisfy the equations~\eqref{eq:d-su3} and hence define the required
invariant forms on~\( \Lambda^2_- \).  

To determine whether a particular form~\( \phi \) given
by~\eqref{eq:su3-phi} is smooth on~\( \Lambda^2_- \), consider the
pull-back of~\( \phi \) to \( P\times S^2H \).  There smoothness reduces to
a question of smooth forms on \( S^2H = \mathbb R^3 \).  Writing these
forms in terms of \( dx_1 \), \( dx_2 \) and \( dx_3 \) one now applies the
results of Glaeser~\cite{Glaeser:compositions}, see also
\cite{Kazdan-W:boundary}, to determine the conditions for the coefficients
to be smooth.  Once \( g \)~is smooth and \( h^*\phi_t=\phi_{-t} \) one
finds that there are no extra conditions.

\providecommand{\bysame}{\leavevmode\hbox to3em{\hrulefill}\thinspace}

\end{document}